# Analytical Solution for Stochastic Unit Commitment Considering Wind Power Uncertainty with Gaussian Mixture Model

Yue Yang, Wenchuan Wu, *Senior Member, IEEE*, Bin Wang, Mingjie Li

*Abstract*—To capture the stochastic characteristics of renewable energy generation output, the chance-constrained unit commitment (CCUC) model is widely used. Conventionally, analytical solution for CCUC is usually based on simplified probability assumption or neglecting some operational constraints, otherwise scenario-based methods are used to approximate probability with heavy computation burden. In this paper, Gaussian mixture model (GMM) is employed to characterize the correlation between wind farms and probability distribution of their forecast errors. In our model, chance constraints including reserve sufficiency and branch power flow bounds are ensured to be satisfied with predetermined probability. To solve this CCUC problem, we propose a Newton method based procedure to acquire the quantiles and transform chance constraints into deterministic constraints. Therefore, the CCUC model is efficiently solved as a mixed-integer quadratic programming problem. Numerical tests are performed on several systems to illustrate efficiency and scalability of the proposed method.

## NOMENCLATURE

### A. Superscripts and subscripts of variables

| | |
|---|---|
| $(.)^t$ | Variables at time $t$ |

### B. Sets and Vectors

| | |
|---|---|
| $\tilde{\mathbf{e}}$ | Random vector composed of forecast errors of wind power from each wind farm |
| $\boldsymbol{\mu}_i$ | Expectation vector of $i$th Gaussian component in GMM |
| $\boldsymbol{\Sigma}_i$ | Covariance matrix of $i$th Gaussian component in GMM |

### C. Variables

| | |
|---|---|
| $T$ | Number of time intervals of optimization |
| $UC_i, FC_i, RC_i$ | UC cost/fuel cost/reserve cost of generator $i$ |
| $SU_i, SD_i$ | Actual startup/shutdown cost of generator $i$ |
| $su_i, sd_i$ | Cost efficient of startup/shutdown of generator $i$ |
| $v_i$ | Commitment state of generator $i$ |
| $a_i, b_i, c_i$ | Coefficients of the fuel cost function of generator $i$ |
| $P_i$ | Power generation schedule of generator $i$ |
| $D_k$ | Power demand of load $k$ |
| $UR_i, DR_i$ | Upward/downward reserve of generator $i$ |
| $urc_i, drc_i$ | Unit cost of upward/downward reserve of generator $i$ |
| $W_{j,cur}$ | Potential wind curtailment of wind farm $j$ |
| $W_{j,forecast}$ | Expected power output based on forecast data of wind farm $j$ |
| $W_{j,sche}$ | Scheduled wind power of wind farm $j$ |
| $\overline{P}_i, \underline{P}_i$ | Upper/lower bound of power generation of generator $i$ |
| $Rup_i, Rdown_i$ | Ramp up/ramp down rate limit of generator $i$ |
| $UR_i^{max}, DR_i^{max}$ | Maximal upward/downward reserve of generator $i$ |
| $UT_i, DT_i$ | Minimal up and down time of generator $i$ |
| $\tilde{w}_j$ | Actual power output of wind farm $j$ |
| $\alpha_{UR}, \alpha_{DR}$ | Maximal allowable probability of upward and downward reserve insufficiency |
| $UR_{extra}, DR_{extra}$ | Extra required upward/downward reserve |
| $s_i^L, s_j^L, s_k^L$ | Power transfer distribution factors |
| $\overline{P}^L$ | Power capacity of transmission line |
| $\alpha^{L+}, \alpha^{L-}$ | Maximal allowable probability of transmission line overloading (bidirectional) |

## I. INTRODUCTION

### A. Background

**W**ITH the integration of massive renewable energy into power system, conventional operation and control schemes based on deterministic optimization encounter the

Manuscript received XX, 2019. This work was supported in part by the National Key R&D Program of China (2018YFB0904200).

Y.Yang,W. Wu and B.Wang are with the State Key Laboratory of Power Systems, Department of Electrical Engineering, Tsinghua University, Beijing 100084, China (e-mail: wuwench@tsinghua.edu.cn). M. Li is with the National Electric Power Control Center of State Grid Corporation of China.

challenges brought by uncertainties of renewable energy sources, which may significantly undermine the reliability and security of power system.

Various approaches have been proposed to extend the deterministic operational procedure by considering uncertainties. As one of those attempts, chance-constrained optimization has been applied in unit commitment (UC) and economic dispatch (ED) problems successfully. In chance-constrained optimization, power outputs of renewable power sources are regarded as random variables with certain probability distribution, and related constraints are required to be satisfied with probability higher than predetermined confidence level to lower the risk of contingencies including reserve exhaustion and transmission line overloading.

Compared to deterministic optimization, solution of chance-constrained optimization is more difficult due to introduction of chance constraints and random variables. The major challenges are (i) accurate probability model of random variables and (ii) efficient transformation from chance constraints to equivalent deterministic ones. Analytical solution method is usually based on simplified probability assumption such as independent Gaussian distribution, while more sophisticated probability model requires iterative scenario-based method to approximate probability by Monte Carlo sampling.

To overcome the above-mentioned difficulties, a chance-constrained UC model is developed in our paper where renewable generation profiles are characterized accurately with Gaussian mixture model and a novel analytical method is proposed to solve the problem efficiently.

*B. Previous research*

Many studies have focused on the application of chance-constrained optimization in power system, including chance-constrained UC and chance-constrained ED problems. A chance-constrained UC model considering load uncertainty is formulated in reference [1], which is solved with iterative method to acquire approximation of chance constraints. Reference [2] and [3] proposed a chance-constrained two-stage UC model to guarantee utilization of wind power and sample average approximation (SAA) algorithm is employed to solve the problem. Authors in [4] presents a chance-constrained UC problem considering the uncertainties of loads and wind power, which is solved by iterations of approximation and verification.

Besides iterative and sampling-based solution methods, there have also been some investigations on efficient transformation from chance-constrained UC/ED model to equivalent deterministic problem. A chance-constrained stochastic programming formulation of day-ahead scheduling is developed and converted into linear deterministic problem in [5]. Reference [6] proposed a reformulation of chance-constrained optimal power flow as second-order cone programming (SOCP) problem by assuming independence and Gaussianity of wind power fluctuation.

On the other hand, some researchers explored the application of non-Gaussian wind power distribution in chance-constrained UC/ED. Versatile distribution is formulated to model the forecast error of wind power and incorporated into chance-constrained ED model in [7]. Reference [8][9] extends the versatile distribution with mixture model to consider correlation between wind farms. Nevertheless, these non-Gaussian models only considers reserve constraints as chance constraints and ignores transmission line constraints due to difficulties to describe the probability distribution of arbitrary linear combination of wind power outputs.

Recently, Gaussian mixture model (GMM) has been adopted to model uncertainties in power systems, including wind power fluctuations[10]-[12], and exhibits promising performance. It can accurately model an arbitrary probability density function (PDF) and possesses favorable mathematical properties facilitating application in terms of chance-constrained optimization problems. A chance-constrained ED model with GMM is developed and solved by approximating the cumulative density function (CDF) of the GMM with the linear combination of piecewise quartic polynomials in [13]. However, this solution has certain flaws: (i) the uniform analytical form of CDF is difficult to obtain because of the piecewise characteristic of each polynomial component; this obstructs a direct solution of the inverse CDF; and, (ii) the existence of multiple solutions to the quartic equation means that an extra validation step is required to obtain a reasonable solution.

*C. Contributions*

In this paper, a stochastic UC model is developed where GMM is introduced to formulate a joint probability distribution for wind farms' output. An analytical solution to acquire the equivalent deterministic formulation of original chance-constrained UC model is proposed. The main contributions of this paper are summarized as follows:

1) We apply GMM to characterize uncertainties of wind power in stochastic UC model. To the best of our knowledge, there is no existing research on the application of GMM in stochastic UC. Compared with Gaussian distribution and other specific probability distributions, GMM is more versatile and can accurately model diverse probability characteristics of wind power fluctuations and correlation between multiple wind farms. Adoption of GMM also facilitate chance-constrained formulation for both reserve constraints and transmission line constraints in our model and subsequent analytical conversion.

2) A novel analytical solution of stochastic UC model is proposed in this paper. By exploiting affine invariance of GMM, the Newton method is applied to efficiently obtain quantiles of one-dimensional GMM and convert all chance constraints to equivalent deterministic linear constraints. The chance-constrained UC formulation is transformed into a mixed-integer quadratic programming (MIQP) problem and solved directly afterwards. Compared to scenario-based solutions of stochastic UC, the proposed method requires no sampling procedure and demonstrates high computational efficiency in numerical tests. Thus, it is applicable for large-scale power systems with renewable energy generation.

3) Potential wind curtailment is considered in our stochastic UC model. In the occasions of insufficient reserve or

transmission capacity, curtailment of wind power is necessary for system security. Thus, introduction of potential wind curtailment ensures the feasibility of UC scheduling. It also provides useful reference for succeeding dispatch and control procedure with finer time scales (as shown in Fig.1) where actual wind curtailment is planned and executed. It is noteworthy that the feasibility of CCUC is seldom discussed in the published analytical solutions, which undermines their applicability.

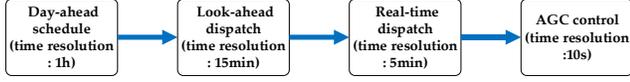

The remainder of this paper is arranged as follows. The chance-constrained formulation of stochastic UC model is introduced in Section II. The solution procedure is discussed in Section III. Numerical tests and results are demonstrated to illustrate the performance of our model in Section IV. Conclusions are drawn in Section V.

## II. MATHEMATICAL MODEL FORMULATION

### A. Objective function

The objective of day-ahead UC optimization is to acquire the optimal UC schedule and generation plan for thermal generators to minimize the overall operation cost which consists of (i) UC cost, (ii) fuel cost, (iii) reserve cost and (iv) potential curtailment penalty of wind power.

$$\min \sum_{t=1}^{T}\left(\sum_{i=1}^{N_G}\left(UC_i^t + FC_i^t + RC_i^t\right) + \sum_{j=1}^{N_W} CP_j^t\right) \quad (1)$$

#### 1) UC cost

The UC cost is composed of startup cost and shutdown cost, which are both dependent on adjacent commitment state of generation units.

$$UC_i^t = SU_i^t + SD_i^t \quad (2)$$
$$SU_i^t = \max\left\{su_i\left(v_i^t - v_i^{t-1}\right), 0\right\} \quad (3)$$
$$SD_i^t = \max\left\{sd_i\left(v_i^{t-1} - v_i^t\right), 0\right\} \quad (4)$$

Here $SU_i^t$ and $SD_i^t$ represent startup cost and shutdown cost respectively and $v_i^t$ is the bool variable that represents the commitment state of generator $i$ at time $t$.

#### 2) Fuel cost

$$FC_i^t = a_i\left(P_i^t\right)^2 + b_i P_i^t + c_i v_i^t \quad (5)$$

When a generator is on ($v_i^t = 1$), fuel cost is described with a quadratic function of its power output where $a_i, b_i, c_i$ denotes coefficients of the fuel cost function. When a generator is off ($v_i^t = 0, P_i^t = 0$), fuel cost is zero.

#### 3) Reserve cost

$$RC_i^t = urc_i UR_i^t + drc_i DR_i^t \quad (6)$$

Reserve cost is composed of cost from upward reserve and downward reserve and both parts are proportional to the corresponding reserve capacity.

#### 4) Potential curtailment penalty of wind power

$$CP_j^t = k_{cur}\left(W_{j,cur}^t\right)^2 \quad (7)$$

Potential curtailment penalty is assumed to be proportional to square of potential curtailed wind power $W_{j,cur}^t$ in order to promote utilization of wind power.

### B. System constraints

#### 1) Power balance constraints

$$\sum_{i=1}^{N_G} P_i^t + \sum_{j=1}^{N_W} W_{j,sche}^t = \sum_{k=1}^{N_D} D_k^t \quad (8)$$

The sum of generation including thermal units and wind farms should equal to sum of power load demand. Here $W_{j,sche}^t$ means the day-ahead scheduled power generation of wind farm $j$ at time $t$.

#### 2) Generator constraints

$$P_i^t + UR_i^t \leq v_i^t \overline{P}_i \quad (9)$$
$$v_i^t \underline{P}_i \leq P_i^t - DR_i^t \quad (10)$$
$$P_i^t - P_i^{t-1} \leq Rup_i + \left(2 - v_i^{t-1} - v_i^t\right)M \quad (11)$$
$$P_i^{t-1} - P_i^t \leq Rdown_i + \left(2 - v_i^{t-1} - v_i^t\right)M \quad (12)$$
$$0 \leq UR_i^t \leq UR_i^{\max} \quad (13)$$
$$0 \leq DR_i^t \leq DR_i^{\max} \quad (14)$$

Equation (9) and (10) are the power output constraints of generators considering upward and downward reserve capacity. Equation (11) and (12) are ramping constraints that are effective if and only if the generator is on at current and previous period with the help of big-M relaxation. Equation (13) and (14) restrict the reserve capacity of individual generator.

It is noteworthy that constraints above are still valid when a generator is off. If generator $i$ is off at time $t$ ($v_i^t = 0$), we will derive $P_i^t = UR_i^t = DR_i^t = 0$ from (9), (10), (13) and (14).

#### 3) Wind power constraints

$$W_{j,sche}^t = W_{j,forecast}^t - W_{j,cur}^t \quad (15)$$
$$\tilde{W}_j^t = W_{j,forecast}^t + \tilde{e}_j^t \quad (16)$$
$$0 \leq W_{j,cur}^t \leq W_{j,forecast}^t \quad (17)$$

Equation (15) gives the relationship between scheduled wind power and its expected forecast value $W_{j,forecast}^t$. In equation (16), $\tilde{W}_j^t$ is the actual power output of wind farm $j$ at time $t$ and is composed of two parts, where $\tilde{e}_j^t$ denotes forecast error and is regarded as random variable. Potential curtailed wind power $W_{j,cur}^t$ is constrained in (17).

#### 4) Minimum up and down time constraints

$$\sum_{k=t}^{t+UT_i-1} v_i^k \geq UT_i\left(v_i^t - v_i^{t-1}\right), \forall t = 1,2,...,T - UT_i + 1 \quad (18)$$

$$\sum_{k=t}^{T}\left[v_i^k - \left(v_i^t - v_i^{t-1}\right)\right] \geq 0, \forall t = T - UT_i + 2,...,T \quad (19)$$

$$\sum_{k=t}^{t+DT_i-1}\left(1-v_i^k\right)\geq DT_i\left(v_i^{t-1}-v_i^t\right), \forall t=1,2,...,T-DT_i+1 \quad (20)$$

$$\sum_{k=t}^{T}\left[1-v_i^k-\left(v_i^{t-1}-v_i^t\right)\right]\geq 0, \forall t=T-DT_i+2,...,T \quad (21)$$

Considering the physical constraints, thermal generators must remain on/off state for several consecutive time periods at least after startup/shutdown. Here $UT_i$ and $DT_i$ denote the minimal required up and down time periods of generator $i$. The mathematical formulation of minimum up and down time constraint in (18)-(21) is based on results from [14] and requires no extra auxiliary variables.

*5) System reserve constraints*

$$\Pr\left(\sum_{i=1}^{N_G}UR_i^t\geq\sum_{j=1}^{N_W}W_{j,sche}^t-\sum_{j=1}^{N_W}\left(\tilde{W}_j^t-W_{j,cur}^t\right)+UR_{extra}\right)\geq 1-\alpha_{UR} \quad (22)$$

$$\Pr\left(\sum_{i=1}^{N_G}DR_i^t\geq-\sum_{j=1}^{N_W}W_{j,sche}^t+\sum_{j=1}^{N_W}\left(\tilde{W}_j^t-W_{j,cur}^t\right)+DR_{extra}\right)\geq 1-\alpha_{DR} \quad (23)$$

When actual wind power output deviates from its scheduled value, there should be enough reserve capacity to maintain system-wide power balance. Due to the randomness of wind power, sufficiency of upward and downward reserve capacity must be guaranteed with high probability, as formulated in (22) and (23).

Here $\alpha_{UR}$ and $\alpha_{DR}$ denote the maximal allowable probability of upward and downward reserve insufficiency. $UR_{extra}$ and $DR_{extra}$ represent the extra required upward/downward reserve margin considering other factors besides wind power fluctuation such as generator outage and load variation.

*6) Transmission line constraints*

$$\Pr\left(\sum_{i=1}^{N_G}s_i^L P_i^t+\sum_{j=1}^{N_W}s_j^L\left(\tilde{W}_j^t-W_{j,cur}^t\right)+\sum_{k=1}^{N_D}s_k^L D_k^t\leq\bar{P}^L\right)\geq 1-\alpha^{L+} \quad (24)$$

$$\Pr\left(\sum_{i=1}^{N_G}s_i^L P_i^t+\sum_{j=1}^{N_W}s_j^L\left(\tilde{W}_j^t-W_{j,cur}^t\right)+\sum_{k=1}^{N_D}s_k^L D_k^t\geq-\bar{P}^L\right)\geq 1-\alpha^{L-} \quad (25)$$

The power flow over transmission line is also stochastic because of uncertainties from wind power. Similar to system reserve constraints, the transmission security constraints are formulated as chance constraints in (24) and (25) to ensure the probability that bidirectional line flow does not exceed its line capacity is bigger than predefined confidence level. Here $s_i^L$, $s_j^L$ and $s_k^L$ denote power transfer distribution factors and $\bar{P}^L$ is the power capacity of transmission line.

## III. SOLUTION PROCEDURE

### A. Modelling of wind power uncertainty

Probability model of wind power uncertainties is vital for stochastic UC scheduling and various distributions have been discussed and applied in previous literatures. Normal distribution is widely used among them to construct probability model of wind power uncertainties, including independent normal distribution for each wind farm in [15][16][17] and multivariate normal distribution as joint distribution for power output of all wind farms in [1][2][3]. Meanwhile, more sophisticated probabilistic modelling of wind power generation based on empirical data analysis such as Copula model is reported in [18][19] but its complexity inhibits the further application in chance-constrained optimization models.

In this paper, Gaussian Mixture Model is employed to characterize uncertainties of wind power. As shown in (26), we assume the joint probabilistic distribution of forecast errors of wind power generation at any moment $t$ denoted by $\tilde{\mathbf{e}}^t$ conforms to a GMM. The PDF of $\tilde{\mathbf{e}}^t$ is convex combination of multiple multivariate Gaussian distributions with corresponding weight coefficient $\alpha_i$, expectation $\boldsymbol{\mu}_i$ and covariance matrix $\boldsymbol{\Sigma}_i$ and each Gaussian distribution is called a component of GMM.

$$PDF_{\tilde{\mathbf{e}}^t}(\mathbf{e})=\sum_{i=1}^{n}\alpha_i N(\mathbf{e},\boldsymbol{\mu}_i,\boldsymbol{\Sigma}_i) \quad (26)$$

$$\tilde{\mathbf{e}}^t=\left[\tilde{e}_1^t,\tilde{e}_2^t,...,\tilde{e}_{N_W}^t\right] \quad (27)$$

$$\sum_{i=1}^{n}\alpha_i=1,\alpha_i>0 \quad (28)$$

The adoption of GMM brings two crucial advantages compared against other probabilistic modelling approach of wind power:

1) The GMM can be used to approach arbitrary probability distribution by adjusting its parameters including number of components, weight coefficients, means and covariance matrixes of each component. Thus, GMM is superior to independent normal distributions and multivariate normal distribution model because it is more flexible and can accurately describe the non-Gaussian characteristics of wind power forecast errors such as asymmetric probabilistic densities. As a multidimensional joint probability distribution, spatial correlation between wind farms are also addressed in GMM.
2) Like multivariate Gaussian distribution, GMM retains the affine invariance, which means the affine transformation of a random variable conforming to GMM still conforms to GMM. The affine invariance property can be exploited to significantly simplify solution of chance-constrained optimization and will be discussed in following derivation. Compared with versatile distribution [7][8] and Copula model, probability distribution of linear combination of random wind power modelled with GMM can be obtained analytically without time-consuming convolution.

### B. Deterministic transformation of chance constraints

By substituting $\tilde{W}_j^t$ from (16), chance constraints in (22)(23)(24)(25) can be transformed to equivalent linear constraints with notation of quantile of the linear combination of random wind power forecast errors. For example, if $\sum_{j=1}^{N_W}\tilde{e}_j^t\left(\mathbf{1}^T\tilde{\mathbf{e}}^t\right)$ is regarded as a new random variable and its $q$-quantile is denoted as $Quant\left(q|\mathbf{1}^T\tilde{\mathbf{e}}^t\right)$, then (22)(23) can be rewritten in following

form:

$$-\sum_{i=1}^{N_G} UR_i^t + UR_{extra} \leq Quant\left(\alpha_{UR} \middle| \mathbf{1^T \tilde{e}}^t \right) \quad (29)$$

$$\sum_{i=1}^{N_G} DR_i^t - DR_{extra} \geq Quant\left(1-\alpha_{DR} \middle| \mathbf{1^T \tilde{e}}^t \right) \quad (30)$$

Similarly, (24)(25) are also rewritten as follows by treating $\sum_{j=1}^{N_W} s_j^L \tilde{e}_j^t \left(\mathbf{s^T \tilde{e}}^t\right)$ as a single random variable:

$$\bar{P}^L - \sum_{i=1}^{N_G} s_i^L P_i^t - \sum_{j=1}^{N_W} s_j^L W_{j,sche}^t - \sum_{k=1}^{N_D} s_k^L D_k^t \geq Quant\left(1-\alpha^{L+} \middle| \mathbf{s^T \tilde{e}}^t \right) \quad (31)$$

$$-\bar{P}^L - \sum_{i=1}^{N_G} s_i^L P_i^t - \sum_{j=1}^{N_W} s_j^L W_{j,sche}^t - \sum_{k=1}^{N_D} s_k^L D_k^t \leq Quant\left(\alpha^{L-} \middle| \mathbf{s^T \tilde{e}}^t \right) \quad (32)$$

*C. Affine invariance and efficient solution of GMM quantiles*

The equivalent form of chance constraints will become deterministic linear constraints after quantiles of $\mathbf{1^T \tilde{e}}^t$ and $\mathbf{s^T \tilde{e}}^t$ are obtained. However, there is no general closed form for the probability distribution of a one-dimensional random variable $\mathbf{c^T \xi}$ where $\mathbf{\xi}$ is a multidimensional random variable and $\mathbf{c}$ is a constant vector even if the probability distribution of $\mathbf{\xi}$ is known and analytical. For example, PDF of $x+y$ where $x$ and $y$ are two independent random variables are described with convolution and computing value of its PDF at arbitrary point requires time-consuming numerical integration.

To obtain the quantiles mentioned above, affine invariance of GMM is introduced and the proof is provided in [20] based on characteristic function.

*Affine invariance of GMM*: If a $n$ dimensional random variable $\mathbf{\xi}$ conforms to GMM and $\mathbf{A}$ is a $n \times m$ dimensional constant matrix, then $\mathbf{A}^T \mathbf{\xi}$ also conforms to GMM as a $m$ dimensional random variable. The PDF of $\mathbf{\xi}$ and $\mathbf{A}^T \mathbf{\xi}$ are given as:

$$PDF_{\mathbf{\xi}}\left(\mathbf{\theta}\right) = \sum_{i=1}^n \alpha_i N\left(\mathbf{\theta}, \mathbf{\mu_i}, \mathbf{\Sigma_i}\right) \quad (33)$$

$$PDF_{\mathbf{A}^T\mathbf{\xi}}\left(\mathbf{\eta}\right) = \sum_{i=1}^n \alpha_i N\left(\mathbf{\eta}, \mathbf{A^T \mu_i}, \mathbf{A^T \Sigma_i A}\right) \quad (34)$$

Affine invariance of GMM implies that $\mathbf{1^T \tilde{e}}^t$ and $\mathbf{s^T \tilde{e}}^t$ both conforms to one dimensional GMM and subsequent derivation will use $\mathbf{s^T \tilde{e}}^t$ as an example. According to (26), PDF and CDF of $\mathbf{s^T \tilde{e}}^t$ are provided as (35) and (36). Here $N(x, \mu_i, \sigma_i^2)$ is the PDF of univariate Gaussian distribution. $\Phi(x)$ is the CDF of the standard Gaussian distribution at $x$.

$$PDF_{\mathbf{s^T \tilde{e}}^t}(x) = \sum_{i=1}^n \alpha_i N\left(x, \mathbf{s^T \mu_i}, \mathbf{s^T \Sigma_i s}\right) \quad (35)$$

$$CDF_{\mathbf{s^T \tilde{e}}^t}(x) = \sum_{i=1}^n \alpha_i \Phi\left(\frac{x - \mathbf{s^T \mu_i}}{\sqrt{\mathbf{s^T \Sigma_i s}}}\right) \quad (36)$$

The $q$-quantile of $\mathbf{s^T \tilde{e}}^t$ is the root of the following univariate nonlinear equation.

$$F(y) = CDF_{\mathbf{s^T \tilde{e}}^t}(y) - q = 0 \quad (37)$$

$F(y)$ increases monotonically and can be efficiently computed when $y$ is given, because $\Phi(y)$ can be obtained directly from the error function. Efficient implementation of the error function has been extensively investigated [21][22] and is available in several software libraries [23][24]. Furthermore, the derivative of $F(y)$ is an elementary function with respect to $y$.

$$F'(y) = PDF_{\mathbf{s^T \tilde{e}}^t}(y) = \sum_{i=1}^n \frac{\alpha_i}{\sqrt{2\pi}\sqrt{\mathbf{s^T \Sigma_i s}}} \exp\left(-\frac{1}{2}\left(\frac{y - \mathbf{s^T \mu_i}}{\sqrt{\mathbf{s^T \Sigma_i s}}}\right)^2\right) \quad (38)$$

Therefore, equation (37) can be iteratively solved using the Newton method.

**Step 1) Initialization**

$$y_0 = y_{init}, k = 0 \quad (39)$$

**Step 2) Iteration**

$$y_{k+1} = y_k - \frac{F(y_k)}{F'(y_k)} \quad (40)$$

**Step 3) Convergence criteria**

If $|F(y_{k+1})| \leq \varepsilon$, the iteration can be terminated and $y_{k+1}$ is the root of equation (37). Otherwise, let $k = k+1$ and return to Step 2 to continue iterating.

To accelerate iteration, the initial value is selected heuristically given the observation that the probability of the quantile in (29)(30)(31)(32) is either close to 1 or close to 0, because allowable reserve exhaustion probabilities $\alpha_{UR}, \alpha_{DR}$ and line overloading probability $\alpha^{L+}, \alpha^{L-}$ are rather small (such as 0.05 or 0.1) in practice. Thus, the initial value of iteration is determined by the probability of quantile $q$:

$$y_{init} = \begin{cases} \max_{i=1...n} \mathbf{s^T \mu_i} & q \geq 0.9 \\ \min_{i=1...n} \mathbf{s^T \mu_i} & q \leq 0.1 \end{cases} \quad (41)$$

**Remark**: The computational performance of Newton method is influenced by the evaluation time of $F(y)$ and its derivative, which are both positively correlated with number of Gaussian components in GMM. As number of wind farms increases, more Gaussian components are essential for accurate probability model and will have negative impact on performance. To accelerate computation, solution of GMM quantiles can be easily parallelized due to its independence between each individual constraint.

*D. Solution framework*

The general solution procedure is shown in flow chart in Fig 2. The construction procedure of GMM to characterize wind power uncertainties will be introduced in IV.A. Chance constraints are turned into deterministic linear constraints as shown in (29)(30)(31)(32) based on quantiles of 1-d GMM obtained via the efficient calculation method in III.C. Afterwards, CCUC model can be solved by existing MIQP solvers

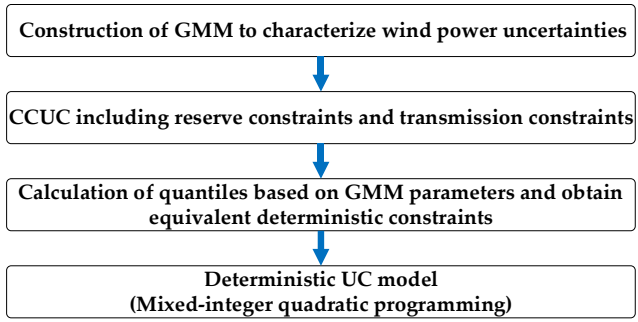

IV. NUMERICAL TESTS

The proposed CCUC formulation with GMM and solution procedure is tested on modified IEEE 24-bus and 118-bus test systems in this paper. All numerical tests were carried out on a laptop with an Intel Core i7-8550U CPU and 16 GB RAM. The calculation of GMM quantiles was implemented in C++ and the deterministic MIQP model was solved by CPLEX (ver. 12.8).

A. *Modified IEEE 24-bus test system*

The modified IEEE 24-bus test system contains 32 thermal generation units and 3 wind farms located at bus 2, 6 and 7 respectively. The total capacity of wind farms is 550MW. The day-ahead unit commitment with 1-hour time resolution is considered in our CCUC model. The maximal probabilities of upward/downward reserve shortage $\alpha_{UR}, \alpha_{DR}$ and transmission line overloading $\alpha^{L+}, \alpha^{L-}$ are all set to 0.02 as default value.

*1) GMM construction for wind power distribution*

In this paper, the construction of GMM to characterize wind power uncertainties is composed of two steps.

In the first step, Nataf transformation is used to generate correlated data samples of wind power forecast errors where forecast error of each wind farm conforms to given marginal probability distribution and predefined correlation coefficients between wind farms. Compared with joint probability distribution, the marginal distributions and correlation coefficients are more accessible in engineering practice. More details of Nataf transformation are discussed in reference [25].

In the second step, Expectation Maximization (EM) algorithm is employed to estimate parameters of GMM based on data samples generated in the first step. We use the implementation of multi-threaded EM in Armadillo C++ library[26] to fit GMM in this paper.

The marginal distribution of forecast error of individual wind farm and correlation coefficients can be estimated with historical wind power profile. In this paper, marginal distribution is modelled by probability histogram. For each time interval, 100000 samples are generated and the number of Gaussian components in EM algorithm is set as 10.

To illustrate the accuracy of GMM compared with normal distribution, the original marginal PDF of forecast errors and marginal PDF obtained from the fitted multivariate normal distribution and GMM are illustrated in Fig. 2. Left skewed, right skewed and bimodal distributions are chosen to address the non-Gaussianity of forecast error distribution. GMM fits these non-Gaussian distributions more precisely in contrast with normal distribution.

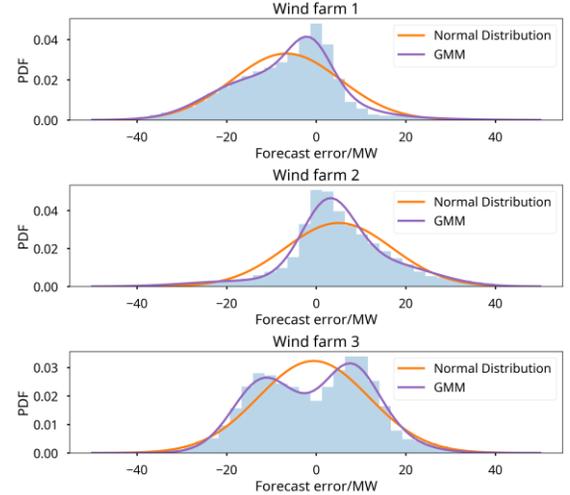

Fig. 3 Original marginal PDF and marginal PDF obtained from normal distribution and GMM for forecast errors in the first time interval

*2) Efficiency and effectiveness of CCUC solution*

The computation burden introduced by chance constraints is evaluated during our optimization procedure by collecting the overall time spent on calculating quantiles of 1-d GMMs to obtain the deterministic equivalence of chance constraints. In our CCUC model for the 24-bus test system, there are 1632 chance constraints and the total time spent on transforming chance constraints is **5.617ms**, namely **3.4μs** per constraint. Meanwhile, the solution time of deterministic MIQP after transformation is beyond **30s** in order to obtain a solution with acceptable relative MIP gap less than 0.01. Thus, quantiles of 1-d GMM are calculated efficiently via our proposed Newton iteration procedure and the extra computation burden of CCUC compared with deterministic UC is relatively ignorable.

To ensure the reserve and transmission line chance constraints are satisfied statistically, Monte Carlo simulation is conducted to validate the commitment schedule. The violation probability of each chance constraint, namely the potential security risk including reserve insufficiency and transmission line overloading is estimated based on randomly generated wind samples.

To illustrate the advantage of GMM over normal distribution, commitment schedules of both probabilistic models are validated through Monte Carlo simulation. The estimated probability of upward reserve insufficiency and overloading of branch from bus 0 to bus 1 is shown in Fig. 4 and Fig. 5. Under UC schedule of GMM, the estimated violation probabilities of chance constraints are restricted approximately under predefined security risk level namely 0.02, whereas UC schedule of normal distribution model fails to satisfy chance constraints in some periods due to its inaccurate probabilistic modelling of forecast errors. Therefore, solution of CCUC with GMM ensures the satisfaction of chance constraints and system security with high probability under wind power uncertainties.

**Remark**: It is noteworthy the CCUC model ignoring potential wind curtailment may encounter infeasible problems when 100% utilization of wind power is restrained by reserve and/or transmission constraints. In this test case, potential wind curtailment occurs during 0:00-2:00 and 6:00-16:00 due to insufficient transmission capacity of branch from bus 7 to bus 9.

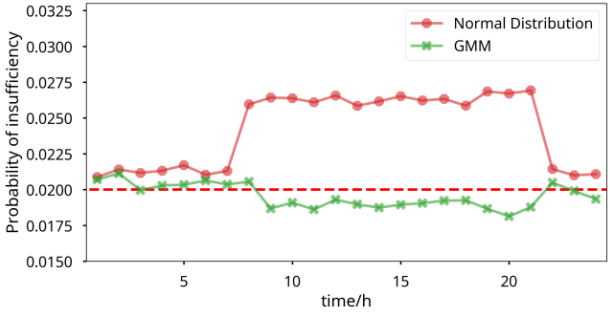

Fig. 4 Probabilities of upward reserve insufficiency estimated by Monte Carlo simulation (24-bus system)

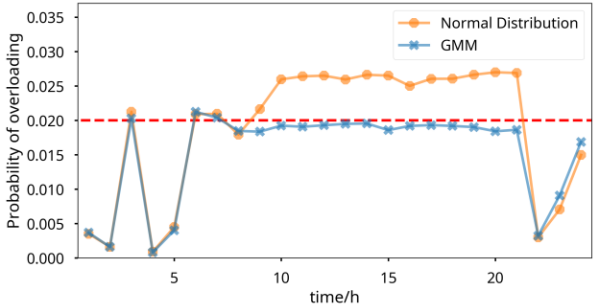

Fig. 5 Probabilities of overloading (branch from bus 0 to 1) estimated by Monte Carlo simulation (24-bus system)

*3) Impact of transmission line constraints*

Transmission capacity is one of the major factors that have remarkable influence on utilization of wind power and optimal UC schedule. Thus, transmission line constraints in form of chance constraints are considered in our model. To illustrate the impacts of chance-constrained transmission limits, Monte Carlo simulation is used to estimate the overloading probability of selected transmission lines under UC schedule with/without consideration of transmission constraints and the result is shown Fig. 6. It is obvious that overloading risk of transmission lines are always kept under predefined level by considering transmission line constraints and will increase significantly if the limits of transmission capacity are ignored in stochastic UC model.

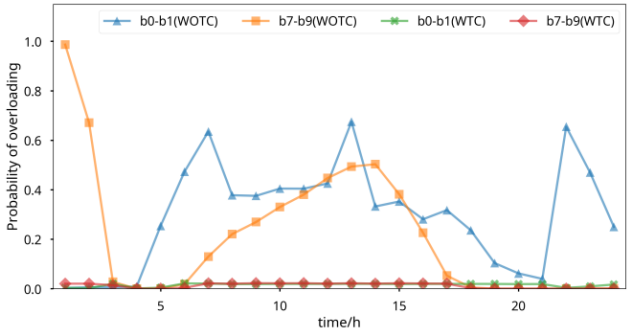

Fig. 6 Probabilities of overloading of selected branches with/without transmission constraints (abbreviated as WTC and WOTC in figure)

*4) Impact of correlation between wind farms*

The uncertainties of different wind farms are correlated due to various factors including spatial climate connections. The correlation structure is described in jointly probability distribution namely GMM of forecast errors and considered in CCUC optimization model. To observe the influence of correlation on optimal UC scheduling, we assume the correlation matrix of 3 wind farms has following form:

$$\begin{bmatrix} 1 & r & r \\ r & 1 & 0 \\ r & 0 & 1 \end{bmatrix} \quad (42)$$

Fig. 7 shows that the overall cost grows monotonically when the correlation coefficient $r$ increases from -0.4 to 0.4. The result can be explained as follows: the forecast errors from different wind farms are offset or augmented by each other when they are negatively or positively correlated respectively. Thus, more reserve from thermal generators and margin of transmission capacity is required to ensure security when correlation coefficient raises, which leads to more conservative UC schedule and higher operational cost.

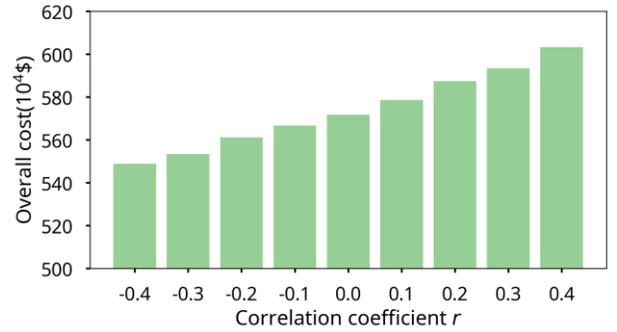

Fig. 7 Overall cost under different correlation coefficients

*B. Modified IEEE 118-bus test system*

There are 54 thermal generation units and 6 wind farms in the modified IEEE 118-bus test system. The total capacity of wind farms is 1200MW.

The entire time spent on computation of GMM quantiles to transform 8016 chance constraints into deterministic form is **108.847ms** namely **13.6μs** per constraint. The time to acquire a solution to MIQP with MIP gap less than 0.01 is beyond **1min**. Compared to the 24-bus system, the average time of each constraint increases as the number of Gaussian components in GMM fitting changes from 10 to 30. The computation burden relative to subsequent MIQP solution is still neglectable, which implies that proposed conversion procedure of chance constraints are applicable to large-scale CCUC problems with numerous chance constraints.

Similar to the 24-bus system test case, Monte Carlo simulation is used to estimate the violation probability of chance constraints. The comparison results between GMM and normal distribution model is shown in Fig. 8 and Fig. 9. The probabilities of upward reserve insufficiency and branch overloading are approximately restricted within predefined value 0.005 during all periods under UC schedule with GMM, which indicates that GMM gives more precise probabilistic description for randomness of wind power than normal distribution and the system security is guaranteed with high probability under proposed stochastic UC schedule.

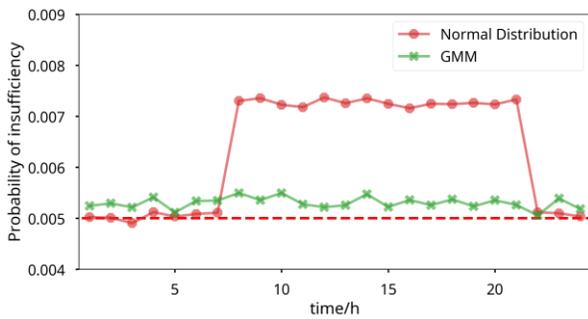

Fig. 8 Probabilities of upward reserve insufficiency estimated by Monte Carlo simulation (118-bus system)

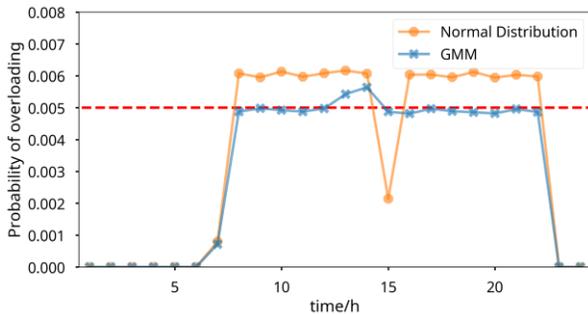

Fig. 9 Probabilities of overloading (branch from bus 7 to 8) estimated by Monte Carlo simulation (118-bus system)

## V. Conclusion

We have proposed a CCUC model where Gaussian mixture model is employed to enable accurate description for diverse probabilistic forecast errors of wind power and correlation structure between wind farms. An efficient solution based on Newton method is proposed to obtain the quantiles of GMMs and transform chance constraints into equivalent deterministic linear constraints. Afterwards, CCUC model is solved as an MIQP problem. Numerical tests are carried out on several test cases. It is shown that GMM fits non-Gaussian correlated distribution of forecast errors precisely and satisfaction of chance constraints under optimal UC schedule is examined by Monte Carlo simulation. Furthermore, our solution method only introduces ignorable computational burden to acquire the equivalent form of CCUC and solve it compared with deterministic UC model. Thus, our proposed modelling and solution to CCUC has potential for large-scale system applications.